\documentclass{crmm-l}

\usepackage{amssymb,amscd,enumerate}

\usepackage{tikz}
\usetikzlibrary{calc,arrows}
\usepackage{diagbox}
\usepackage{pdfsync}
\usepackage{caption}

\usepackage{chngcntr}
\counterwithout{figure}{chapter}

\usepackage[cmtip,all]{xy}

\usepackage[bookmarks, colorlinks=true, pdftitle={Corrigendum to: Elliptic Boundary Value Problems with Fractional Regularity Data}, pdfauthor={Alex Amenta, Pascal Auscher and Moritz Egert}]{hyperref}

\makeatletter
\@namedef{subjclassname@2020}{\textup{2020} Mathematics Subject Classification}
\makeatother

\theoremstyle{definition}

\theoremstyle{remark}

\numberwithin{section}{chapter}
\numberwithin{equation}{chapter}

\newcommand{\mc}{\mathcal}
\newcommand{\mb}{\mathbf}
\newcommand{\nm}[1]{\| #1 \|} 

\newcommand{\ga}{\alpha}

\newcommand{\gd}{\delta}

\newcommand{\gf}{\varphi}

\newcommand{\gh}{\eta}

\newcommand{\gq}{\theta}

\newcommand{\gs}{\sigma}
\newcommand{\gt}{\tau}

\newcommand{\gy}{\psi}

\newcommand{\gY}{\Psi}

\newcommand{\bbR}{\mathbb{R}}

\def\barint_#1^#2{\mathchoice
            {\mathop{\vrule width 4.5pt
height 3 pt depth -2.5pt
                    \kern -7.4pt
\int_{#1}^{#2} \kern 0pt}}%
            {\mathop{\vrule width 5pt height
3 pt depth -2.6pt
                    \kern -6.5pt
\intop_{#1}^{#2} \kern -4pt}\nolimits_{#1}}%
            {\mathop{\vrule width 5pt height
3 pt depth -2.6pt
                    \kern -6pt
\intop_{#1}^{#2} \kern -4pt}\nolimits_{#1}}%
            {\mathop{\vrule width 5pt height
3 pt depth -2.6pt
          \kern -6pt \intop \kern -4pt}\nolimits_{#1}}}
          
\def\bariint_#1^#2{\mathchoice
            {\mathop{\vrule width 8pt
height 3 pt depth -2.5pt
                    \kern -10.5pt
\intop \kern -9pt\int_{#1}^{#2} \kern 0pt}}%
            {\mathop{\vrule width 9pt height
3 pt depth -2.6pt
                    \kern -10.5pt
\intop \kern -10pt\int_{#1}^{#2} \kern -4pt}}%
            {\mathop{\vrule width 9pt height
3 pt depth -2.6pt
                    \kern -10pt
\intop \kern -10pt\int_{#1}^{#2} \kern -4pt}}%
            {\mathop{\vrule width 9pt height
3 pt depth -2.6pt
          \kern -10pt \int_{#1}^{#2} \kern -10pt\intop \kern -4pt}
      }}
          
\newcommand{\dint}{\int \kern -9pt\int}

\makeindex

\allowdisplaybreaks

\frontmatter

\title[Elliptic BVPs with Fractional Regularity Data]{Corrigendum to: \\ Elliptic Boundary Value Problems with Fractional Regularity Data: The First Order Approach}

\author{Alex Amenta}
\address{}
\email{amenta@fastmail.fm}

\author{Pascal Auscher}
\address{Universit\'e Paris-Saclay, CNRS, Laboratoire de Math\'{e}matiques d'Orsay, 91405 Orsay, France}
\email{pascal.auscher@universite.paris.saclay.fr}

\author{Moritz Egert}
\address{Fachbereich Mathematik, Technische Universit\"at Darmstadt, Schlossgartenstr. 7, 64289 Darmstadt, Germany}
\email{egert@mathematik.tu-darmstadt.de}



\begin{document}
\begin{abstract}
This is a corrigendum to the monograph  ``Elliptic Boundary Value Problems with Fractional Regularity Data: The First Order Approach'' by the first two authors. It addresses two major imprecisions in the preliminary material. 
\end{abstract}

\maketitle
\pagenumbering{arabic} 

The preliminary material of the monograph \cite{AA} written by the first two authors contains two major imprecisions that necessitates a number of (in the end harmless) changes throughout the entire text. One is about identification of abstract and concrete Hardy spaces for perturbed Dirac operators, the other one about interpolation of quasi-Banach function spaces. Since these erroneous statements are not unlikely to spread, we provide a detailed corrigendum, including further background and corrected statements for all affected results. All other results remain unchanged.

\section*{A correction on the classification region of Section 6}

In the introduction, Theorem~1.8,  summarising references \cite{AS16} and  \cite{AM15}, is stated  incorrectly. This induces some changes in later statements in the introduction and the definition of the classification region in Chapter~6 should be modified for the results in Chapter~6 and 7 to be fully correct. 

As is mentioned after the statement of Theorem~1.8, there is an interval of exponents that has been introduced in \cite{AS16} and every $p$ in this interval is such that the assumption $\mb{H}^p_{DB} \simeq \mb{H}_D^p$ is satisfied. It was unfortunate that this interval was named  $I_0(\mb{H},DB)$, because this notation is used in Chapter 5 with another definition (after Definition 5.10) as the maximal interval in $(0,\infty)$ for which $\mb{H}^p_{DB} \simeq \mb{H}_D^p$,  causing the confusion. A detailed explanation of the possible difference is given \cite[Chapter 15]{AE}, to which we refer. 

Let us call $p_{+}(DB)$ the upper bound of the interval introduced \cite{AS16}. It is a number in $(2,\infty]$ with a precise meaning. Its definition is used in the proof of uniqueness obtained in  \cite{AM15} and therein the results are (correctly) stated with the constraint $p<p_{+}(DB)$. Consequently, the following should be corrected here.

\begin{quote}
	\begin{itemize}
		
		\item One  should correct  Theorem~1.8 by assuming
		that $p$ in  (i) satisfies $p<p_{+}(DB)$ and that $p'$ in (ii) satisfies $p'<p_{+}(DB^*)$. 
		\item Next, Theorem~1.11 should be modified by assuming $p<p_{+}(DB)$ when $\theta=0$ and $p>(p_{+}(DB^*)'$ when $\theta=-1$. 
		\item The same modification should be applied in Theorems~1.12, case (i), and 1.14,  to the pairs $(\theta_{0},p_{0})$ and $(\theta_{1},p_{1})$. 
	\end{itemize}
\end{quote}
Since all identification results for $\theta(\mb{p}) \in \{0,1\}$ that appear in the later chapters refer back to \cite{AS16, AM15}, also the definition of the classification region $J(\mb{H},DB)$ in Definition~6.10 should be adapted. It should read as follows:

\begin{quote}
	\textbf{Definition 6.10.} We define the \emph{classification region} for $DB$ as follows. For $\mb{X}=\mb{B}$ we let \pagebreak
	\begin{align*}
		J(\mb{B},DB) := \{\mb{p} \in I_{\text{max}} : &\text{[$i(\mb{p}) \leq 2$ and $\mb{p} \in I(\mb{B},DB)$]} \\
		&\text{or [$i(\mb{p}) > 2$ and $\mb{p}^\heartsuit \in I(\mb{B},DB^*)$]}\}
	\end{align*}
	and for $\mb{X}=\mb{H}$ we let
	\begin{align*}
		\qquad \qquad J(\mb{H},DB) \\:= \{\mb{p} &\in I_{\text{max}} : \text{[$-1<\theta(\mb{p})<0$],} \\
		&\qquad \qquad \quad  \text{and  [$i(\mb{p}) \leq 2$ and $\mb{p} \in I(\mb{H},DB)$]} \\
		&\qquad \qquad \quad \text{or [$i(\mb{p}) > 2$ and $\mb{p}^\heartsuit \in I(\mb{H},DB^*)$],}
		\\ 
		&\hphantom{\in I_{\text{max}} : \text{}}\text{ [$\theta(\mb{p})=0]$,} \\
		&\qquad \qquad \quad \text{and  [$i(\mb{p}) \in I_0(\mb{H},DB)$, $i(\mb{p})<p_{+}(DB)$],}
		\\ 
		&\hphantom{\in I_{\text{max}} : \text{}} \text{[$\theta(\mb{p})=-1]$,}\\
		&\qquad \qquad \quad  \text{and [$i(\mb{p}^\heartsuit) \in I_{0}(\mb{H},DB^*)$, $i(\mb{p}^\heartsuit) <p_{+}(DB^*)$]}\}.
	\end{align*}
\end{quote}

We note that there is no modification for the classification region  $J(\mb{B},DB)$ for Besov data, as it is only concerned with regularity $-1<\theta(\mb{p})<0$. The correct classification region $J(\mb{H},DB)$ is smaller than what has originally been stated in the monograph. Consequently, the results come with additional restrictions when $\theta(\mb{p}) \in \{0,1\}$, which, again, were treated in \cite{AM15}, whereas proofs when $0<\theta(\mb{p}) <1$ remain correct as they stand. 

Finally,  let us mention the question whether there exists $p\in [p_{+}(DB),\infty)$  for which $\mb{H}^p_{DB} \simeq \mb{H}_D^p$. In the case where $B$ is a block matrix,\cite[Chapter 15]{AE} establishes that the answer is no.  The general case remains not fully understood at this stage. 

\section*{A correction on interpolation results}

The interpolation results, in particular in Section~2.3, require two main adjustments. To this end, it is instructive to consider the exponent $\mb{p}$ in the $(j,\theta)$-plane and distinguish three regions as in the figure below: The finite range $j(\mb{p}) > 0$, the infinite range $j(\mb{p}) \leq 0$ and the duality range $0 \leq j(\mb{p}) < 1$.

\begin{figure}[ht]\label{fig:exponents_regions}

\begin{center}
\begin{tikzpicture}[scale=2]

	\draw [thick, dashed] (2,1.5) -- (2,-2); 
	\path [fill=lightgray, opacity = 0.4] (0,1.5) -- (2.,1.5) -- (2.,-2) -- (0,-2);
	
	\draw [thick,->] (-0.25,0) -- (3,0); 
	\draw [thick,->] (0,-2) -- (0,1.5); 
	
	\draw [fill=black] (0.75,-0.5) circle [radius = 1pt]; 
	\node [right] at (0.75,-0.5) {$\mb{p}$};

	\draw [fill=black] (0.5,-1.25) circle [radius = 1pt]; 
	\node [right] at (0.5,-1.25) {$\mb{q}$};
	
	\draw [fill=black] (0.625,-0.875) circle [radius = 0.75pt]; 
	\node [right] at (0.625,-0.875) {$[\mb{p},\mb{q}]_{\theta}$};
	
	\draw[dashed] (0.75,-0.5) -- (0.5,-1.25);

	\node [above] at (0,1.5) {$\theta$};
	\node [left] at (0,1) {$\frac{1}{2}$};
	\node [left] at (0,-1) {$-\frac{1}{2}$};
	\node [left] at (0,-2) {$-1$};
	\node [below] at (1.1,0) {$\frac{1}{2}$};
	\node [below] at (2.1,0) {$1$};
	\node [below] at (2.66,0) {$\frac{n+1}{n}$};
	\node [right] at (3,0) {$j$};
\end{tikzpicture}
\end{center}
	\captionsetup{labelformat=empty}
\caption{\textsc{Exponent regions}: The gray duality range described by $0 \leq j(\mb{p}) < 1$ corresponds to the spaces that can be obtained as duals of Banach spaces with finite exponent $\mb{p}$ (that is, $j(\mb{p}) > 0$). It lies in the intersection of finite exponents and infinite exponents with $j(\mb{p}) = 0$. The geometrically interpolating exponent $[\mb{p},\mb{q}]_{\theta}$ captures the interpolation behaviour of function spaces if $j(\mb{p}), j(\mb{q}) \geq 0$.} \label{fig:exponents}
\end{figure}

The first adjustment is that interpolation of tent spaces, $Z$-spaces and homogeneous smoothness spaces only translates to taking convex combinations of exponents in the $(j,\theta)$-plane, if we restrict ourselves to the right half-plane with $j(\mb{p}) \geq 0$. 

For example, real interpolation of a Besov space $ \dot{B}^{\infty,\infty}_{s + \ga}$ with $\dot{B}^{p,p}_{t}$  for finite $p$ always yields a Besov space with finite exponent --- no matter how far to the left the exponent $\mb{p} = (\infty,s;\alpha)$, corresponding to $(-\alpha/n,s)$, lies  in the $(j,\theta)$-plane, compare with the references in Section~2.3. 

The additional assumption $j(\mb{p}) \geq 0$ is a (harmless) restriction for results on tent and Hardy--Sobolev spaces. 

For $Z$-spaces, we see from Definition~2.15 that $Z^{\infty,s;\alpha} = Z^{\infty,s+\alpha;0}$, and the same applies to Besov spaces. Hence, the same space can have different parametrisations in the $(j,\theta)$-plane and for a correct formulation of interpolation results we have to pick the unique one with $j(\mb{p})\geq 0$. In other words, when $\mb{p}, \mb{q}$ are  infinite with $0=j(\mb{p})$ and $\mb{p}\hookrightarrow \mb{q}$, then we have by definition that $Z^{\mb{q}}=Z^{\mb{p}}$. The distinct $Z$-spaces are only those with $j(\mb{p})\geq 0$. The same applies for the Besov scale. 

The second adjustment is connected to a general strategy of proof for interpolation results in three steps: first one checks the finite range, then gets the duality range by duality methods for interpolation functors, and finally one gets all exponents by Wolff reiteration because the finite range and the duality range have a substantial overlap. The issue is that no Wolff reiteration is known for the complex interpolation method of Kalton--Mitrea!

The note~\cite{Egert-Kosmala} contains a proof of the usual Wolff reiteration theorem for four $p$-convex quasi-Banach function spaces, where only one of the `exterior spaces' needs to be separable, provided that all four spaces have the so-called Fatou property. For details, we refer to \cite{Egert-Kosmala}, whereas for the purpose of this corrigendum it suffices to note that tent- and $Z$-spaces are covered if one of the `exterior exponents' is finite, which is the case in the scenario described above. This is how Wolff reiteration for complex interpolation spaces has to be understood in the context of Section~2.3.

We come to the concrete changes in the text that follow from the above two main adjustments and their consequences. All other results in the monograph remain true as stated with the sole exception of Theorem~7.8, for which we describe the modification further below.

\section*{Changes in Section~2.3}

The corrected statement and proof for real interpolation is as follows. We need to assume $j(\mb{p}), j(\mb{q}) \geq 0$, also in order to justify duality for real interpolation.

\begin{quote}
	\textbf{Theorem 2.30 (Real interpolation of tent spaces: full range).} 
	\textit{Suppose that $\mb{p}$ and $\mb{q}$ are exponents with $\gq(\mb{p}) \neq \gq(\mb{q})$ and $j(\mb{p}), j(\mb{q}) \geq 0$, and $0 < \gh < 1$.
		Then
		\begin{equation*}
			(T^\mb{p}, T^\mb{q})_{\gh,p_\gh} = Z^{[\mb{p},\mb{q}]_\gh},
		\end{equation*}
		where $p_\gh = i([\mb{p},\mb{q}]_\gh)$.}
	
	\begin{proof}
		For finite exponents this is Theorem~2.18.
		Let now $1 < i(\mb{p}), i(\mb{q}) \leq \infty$. Since $j(\mb{p}), j(\mb{q}) \geq 0$, we have that $T^{\mb{p}^\prime}, T^{\mb{q}^\prime}$ are Banach spaces. Hence, in this case the claim follows by writing
		\begin{equation*}
			(T^{\mb{p}}, T^{\mb{q}})_{\gh,p_\gh} = ((T^{\mb{p}^\prime})^\prime, (T^{\mb{q}^\prime})^\prime)_{\gh,p_\gh} = (T^{\mb{p}^\prime}, T^{\mb{q}^\prime})_{\gh,p_\gh^\prime}^\prime
		\end{equation*}
		via the duality theorem for real interpolation [22, Theorem 3.7.1], using that $T^{\mb{p}^\prime} \cap T^{\mb{q}^\prime}$ is dense in both $T^{\mb{p}^\prime}$ and $T^{\mb{q}^\prime}$, and then noting that
		\begin{equation*}
			p_\gh^\prime = i([\mb{p},\mb{q}]_\gh)^\prime = i([\mb{p}^\prime, \mb{q}^\prime]_\gh).
		\end{equation*}
		The full result follows by Wolff reiteration~[76, Theorem 1] which, for the real method holds in the quasi-Banach setting.
	\end{proof}
\end{quote}

In the subsequent result on real interpolation, the same strategy has to be used. For finite exponents one uses reiteration (but the one allowing quasi-Banach spaces in  [22, Theorem 3.11.4]). Then duality and Wolff reiteration complete the argument as before. The assumption $j(\mb{p}), j(\mb{q}) \geq 0$ is needed again, however, as in Theorem~2.30, they can both be infinite:

\begin{quote}
	\textbf{Proposition 2.31 (Real interpolation of $Z$-spaces).}
	\textit{Let $\mb{p}$ and $\mb{q}$ be exponents with  $\gq(\mb{p}) \neq \gq(\mb{q})$ and $j(\mb{p}), j(\mb{q}) \geq 0$, and let $\gh \in (0,1)$.
		Then
		\begin{equation*}
			(Z^\mb{p}, Z^\mb{q})_{\gh,p_\gh} = Z^{[\mb{p},\mb{q}]_\gh},
		\end{equation*}
		where $p_\gh = i([\mb{p}, \mb{q}]_\gh)$.}
\end{quote}

Also in Proposition~2.33 we need $j(\mb{p}), j(\mb{q}) \geq 0$. The strategy is once again the same and the case of finite exponents is what is being deferred to Section~2.6. Duality for the complex method on Banach spaces requires one space to be reflexive, which explains why the assumption that one exponent is finite has to be maintained as in the original text.

\begin{quote}
	\textbf{Proposition~2.32 (Complex interpolation of $Z$-spaces).}
	\textit{Let $\mb{p}$ and $\mb{q}$ be exponents with $j(\mb{p}), j(\mb{q}) \geq 0$ and equality for at most one of them, and let $\gh \in (0,1)$.
		Then
		\begin{equation*}
			[Z^\mb{p}, Z^\mb{q}]_{\gh} = Z^{[\mb{p},\mb{q}]_\gh}.
	\end{equation*}}
\end{quote}

The restriction to exponents with $j(\mb{p}) \geq 0$ also requires some more care in the proof of the mixed embedding in Theorem~2.34 because in the proof $\mb{r} = [\mb{p}, \mb{q}]_2$ may not satisfy $j(\mb{r}) \geq 0$ even if $\mb{p}$ and $\mb{q}$ are finite. A three step strategy as above is once again necessary and we give a corrected proof of the unchanged statement. 

\begin{quote}
	\textbf{Theorem 2.34 (Mixed embeddings).}
	\textit{Let $X_0,X_1 \in \{T,Z\}$ and let $\mb{p} \hookrightarrow \mb{q}$ with $\mb{p} \neq \mb{q}$.
		Then we have the embedding}
	\begin{equation*}
		(X_0)^\mb{p} \hookrightarrow (X_1)^\mb{q}.
	\end{equation*}
	\begin{proof}
		We first treat the case that $\mb{p}$ and $\mb{q}$ are finite. When $X_0 = X_1 = T$, this is Theorem~2.11.
		
		Let $\mb{r} = [\mb{p}, \mb{q}]_{1+\epsilon}$ with $\epsilon>0$ small enough to guarantee that $\mb{r}$ is still finite. Then $\mb{p} \hookrightarrow \mb{r}$ and $[\mb{p},\mb{r}]_{(1+\epsilon)^{-1}} = \mb{q}$ (by Lemmas~2.1 and 2.2).
		Then we have embeddings $T^\mb{p} \hookrightarrow T^{\mb{p}}$ (trivially) and $T^\mb{p} \hookrightarrow T^\mb{r}$ (Theorem~2.11).
		Hence
		\begin{equation*}
			T^\mb{p} \hookrightarrow (T^\mb{p}, T^\mb{r})_{(1+\epsilon)^{-1},i([\mb{p},\mb{r}]_{(1+\epsilon)^{-1}})} = Z^{[\mb{p},\mb{r}]_{(1+\epsilon)^{-1}}} = Z^\mb{q}
		\end{equation*}
		by Theorem~2.30, using that $\mb{p} \neq \mb{q}$ and $\mb{p} \hookrightarrow \mb{q}$ imply $\gq(\mb{p}) \neq \gq(\mb{q})$.
		Similarly, putting $\mb{s} = [\mb{p},\mb{q}]_{-1}$, we have $T^\mb{s} \hookrightarrow T^{\mb{q}}$ and $T^{\mb{q}} \hookrightarrow T^{\mb{q}}$, so
		\begin{equation*}
			Z^\mb{p} = (T^\mb{s}, T^\mb{q})_{1/2,i([\mb{s},\mb{q}]_{1/2})} \hookrightarrow T^\mb{q}.
		\end{equation*}
		Finally, putting $\mb{t} = [\mb{p},\mb{q}]_{1/2}$ and using the previous results, we have
		\begin{equation*}
			Z^\mb{p} \hookrightarrow T^\mb{t} \hookrightarrow Z^\mb{q},
		\end{equation*}
		which completes the proof for finite exponents.
		
		The next case is that $\mb{p'}$ and $\mb{q'}$ are both finite. Then the claim follows from the first case by duality and this covers the infinite exponents and those in the duality range, see the figure above.
		
		The remaining case is that $\mb{q}$ is infinite and $\mb{p}$ is finite but not in the duality range. Then the connecting segment contains a finite exponent $\mb{r} =  [\mb{p}, \mb{q}]_{\theta}$ with $\theta \in (0,1)$ in the duality range. Combining the previous two cases yields
		\begin{equation*}
			(X_0)^\mb{p} \hookrightarrow (X_0)^\mb{r}  \hookrightarrow (X_1)^\mb{r}.
		\end{equation*}
		The proof is complete.
	\end{proof}
\end{quote}
\section*{Changes in Section~4.3}

The interpolation results of Section~2.3 directly correspond to the interpolation results for adapted spaces in Theorem~4.28. Hence, this theorem should be read as follows:

\begin{quote}
	\textbf{Theorem 4.28 (Interpolation of completions).}
	Fix $0 < \gh < 1$ and $\gy \in \gY_\infty^\infty$.
	Let $\mb{p}_0$ and $\mb{p}_1$ be exponents, and set $\mb{p}_\gh := [\mb{p}_0,\mb{p}_1]_\gh$.
	\begin{enumerate}[(i)]
		\item
		Suppose $j(\mb{p}_0), j(\mb{p}_1) \geq 0$, with equality for at most one exponent.
		Then we have the identification
		\begin{equation*}
			[\gy \mb{H}_A^{\mb{p}_0}, \gy \mb{H}_A^{\mb{p}_1}]_\gh = \gy\mb{H}_A^{\mb{p}_\gh}.
		\end{equation*}
		
		\item
		Suppose $j(\mb{p}_0), j(\mb{p}_1) \geq 0$, with equality for at most one exponent.
		Then we have the identification
		\begin{equation*}
			[\gy \mb{B}_A^{\mb{p}_0}, \gy \mb{B}_A^{\mb{p}_1}]_\gh = \gy\mb{B}_A^{\mb{p}_\gh}.
		\end{equation*}
		
		\item
		Suppose $j(\mb{p}_0), j(\mb{p}_1) \geq 0$ and that $\gq(\mb{p}_0) \neq \gq(\mb{p}_1)$.
		Then we have the identification
		\begin{equation*}
			(\gy \mb{X}_A^{\mb{p}_0}, \gy \mb{X}_A^{\mb{p}_1})_{\gh,i(\mb{p}_\gh)} = \gy\mb{B}_A^{\mb{p}_\gh}.
		\end{equation*}
	\end{enumerate}
	All of these statements have analogues for spectral subspaces.
\end{quote}

\section*{Changes in Section~6.5}

The interpolation results of Section~2.3 and Section~4.3 have been joined in the proof of Theorem~6.41. Hence, also in this theorem, (ii) needs the additional assumption $j(\mb{p_1}) \geq 0$ and in (iii) it is needed that $j(\mb{p_0}), j(\mb{p_1}) \geq 0$.

\section*{Changes in Section~7.1}

Theorem~7.8 is the only result, where infinite exponents mistakenly appear by interpolation between infinite and finite exponent. In fact, {\v{S}}ne{\u\i}berg's extrapolation theorem [66] only applies to interior points of complex interpolation scales. This result only works for finite exponents. On the other hand, the restriction $i(\mb{p})>1$ in the Besov case is unnecessary because complex interpolation for $Z$-spaces was also obtained in the quasi-Banach range in Proposition~2.32. (Section~2.6 was finalised at the very end of composing the monograph.) Theorem~7.8 and the succeeding remark should therefore be modified as follows.

\begin{quote}
	\textbf{Theorem 7.28 (Extrapolation of well-posedness)}
	\textit{Let $B = \hat{A}$, and suppose $\mb{p} \in J(\mb{X},DB)$ is a finite exponent with $\gq(\mb{p}) \in (-1,0)$.
		Suppose that $(P_\mb{X})_A^\mb{p}$ is well-posed.
		Then there exists a $(j,\gq)$-neighbourhood $O_\mb{p}$ of $\mb{p}$ such that for all $\mb{q} \in O_\mb{p}$ with $\theta(\mb{q}) \leq 0$, $(P_\mb{X})_A^{\mb{q}}$ and $(P_\mb{X})_A^\mb{p}$ are mutually well-posed.}
\end{quote}

\medbreak

\begin{quote}
	\textbf{Remark 7.29}
	\textit{With $\mb{X} = \mb{H}$, a corresponding result is true for $\gq(\mb{p}) \in \{-1,0\}$, but well-posedness is only obtained for $\mb{q}$ near $\mb{p}$ with $\gq(\mb{q}) = \gq(\mb{p})$. The proof uses the same argument.}
\end{quote}

It should be noted that the subsequent results (Corollary~7.10 and Theorem~7.11) are not affected by the modification above since they only use finite exponents anyway.

\section*{Further typos}

\begin{itemize}
	
	\setlength\itemsep{0.8em}
	
	\item Front matter: in the Library of Congress Cataloging-in-Publication Data, the first author's year of birth is off by 10 years. The direction of the error is left as an exercise to the reader
	
	\item Page 4, line 2: ``Holder" should read ``H\"older".
	
	\item In the line before Theorem~2.11, replace ``by the same argument''   with ``by duality, using Theorem~2.7''.
	
	\item  Page 24: In the formula (2.7) for $\nm{f}_{Z^\mb{p}_c}$, it is implicit that the measure on $\bbR^{1+n}_+$ is $\frac{dxdt}t$. 
	
	\item Page 25: In Proposition~2.20, change $\ell^p$ to $\ell^{i(\mb{p})}$. 
	
	\item Page 30, line 5: $\nm{g}_{Z_c^\mb{p}}$ should read $\nm{g}_{Z_c^{\mb{p}'}}$. 
	
	\item Page 31, line before Proposition~2.31: The correct reference is [23, Theorem~3.11.5].
	
	\item Page 36, line -7: Change $D^\alpha f(0)$ to $D^\alpha \hat f(0)$.
	
	\item Page 37, statement of Theorem~2.49: Change $s$ to $\alpha$.
	
	\item Page 56, Lemma~3.17:  The statement should be modified to 
	$\gy \in \gY_{\gs}^{\gt+}$, $\gf \in \gY_{\gt + \gd}^{(\gs - \gd)+}$. In the proof, when $\gs + \gt = 0$ there is uniform boundedness, and Theorem~3.8 applies when $\gs + \gt > 0$.
	
	\item Page 66: In the definition of a completion, replace ``continuous injective'' by ``isometric''.
	
	\item Page 69, Proposition~4.24: For infinite exponents, bounded means continuous for the weak$^*$ topology (obtained by duality). For the density in the proof, one can select a sequence $(F_{k})$ in the proof of Corollary~4.9 with compact support in $t$, so that $f_{k} \in \mc{D}(\gh(A))$. 
	
	\item Page 69, Remark~4.25: There is also a completed version of the density property in Corollary~4.9. This is necessary to perform interpolation.
	
	\item Page 76: In the second display of the proof of Theorem 5.3, one should replace $t^{-1}\,\cdot$ by $t^{-1}|\cdot|$.
	
	\item Pages 85-87, proof of Lemma 5.29: There are several adjustments necessary.
	
	(a) The value of $const.$ in  line 2 and 3 should be $1$. The assumptions that   $\|\chi\|_{\infty}+\tau\|\nabla \chi\|_{\infty} \lesssim 1$ and $\|\chi_{1}\|_{\infty}+\tau\|\nabla \chi_{1}\|_{\infty} \lesssim 1$ are missing. 
	
	(b) In the series of inequalities at the bottom of page 86, on the fourth inequality, one should have used the Cauchy--Schwarz inequality for the integral in $y$ and then argue as suggested. Moreover, $r^n$ should read $\tau^n$ (as on the line after). 
	
	(c) In the series of inequalities at the top of page 87, in the third inequality a factor $\tau^{-2}$ is missing. Moreover, one should delete the subscript $\chi$ and in the fourth inequality, replace $r^{-2\theta}$ by $\tau^{-2\theta-2}$. 
	
	\item Page 95, Proposition~6.1: The Whitney parameter on the left-hand side should be smaller than $c_{0}$ and not equal to $c_{0}$.
	
	\item Page 98, Lemma~6.7: The limit is $t\to \infty$, not $t\to 0$. 
	
	\item Page 101, line -9: The convergence below the display is only weakly in $L_\text{loc}^2(\bbR^{1+n}_+)$ for infinite exponents.
	
	\item Page 117, Corollary~6.38: ${\mb{X}}_{D}^{\mb{p},+}$ should read ${\mb{X}}_{DB}^{\mb{p},+}$.
	
	\item Page 128, Theorem 7.5: The isomorphisms should be understood as topological isomorphisms (hence, bi-continuous) and the open mapping theorem is used implicitly in the proof. In the quasi-Banach range one reference is \cite[Theorem~1.4]{Kalton-Peck-Roberts}. In the infinite exponent range, continuity is for the weak $^*$-topology, which in all appearing cases is a Fr\'echet topology.  
	
	\item Page 129, line 10: In the sentence starting by ``evidently'', one should mention the use of a density argument of three spaces of the considered type in the intersection of two of them, together with the continuity of the inverses proved in Theorem 7.5.

\end{itemize}


\backmatter
\bibliographystyle{amsplain}
\bibliography{amenta-auscher-bib}

\providecommand{\bysame}{\leavevmode\hbox to3em{\hrulefill}\thinspace}
\providecommand{\MR}{\relax\ifhmode\unskip\space\fi MR }
\providecommand{\MRhref}[2]{%
  \href{http://www.ams.org/mathscinet-getitem?mr=#1}{#2}
}
\providecommand{\href}[2]{#2}
\begin{thebibliography}{1}

\bibitem{AA}
A.~Amenta and P.~Auscher, \emph{Elliptic boundary value problems with
  fractional regularity data}, CRM Monograph Series, vol.~37, American
  Mathematical Society, Providence, RI, 2018, The first order approach.
  \MR{3753666}

\bibitem{AE}
P.~Auscher and M.~Egert, \emph{Boundary value problems and {H}ardy spaces for
  elliptic systems with block structure}, Progress in Mathematics, vol. 346,
  Birkh\"{a}user/Springer, Cham, 2023. \MR{4628043}

\bibitem{AM15}
P.~Auscher and M.~Mourgoglou, \emph{Representation and uniqueness for boundary
  value elliptic problems via first order systems}, Rev. Mat. Iberoam.
  \textbf{35} (2019), no.~1, 241--315.

\bibitem{AS16}
P.~Auscher and S.~Stahlhut, \emph{Functional calculus for first order systems
  of {Dirac} type and boundary value problems}, M\'emoires de la Soci\'et\'e
  Math\'ematique de France, vol. 144, Soci\'et\'e Math\'ematique de France,
  Paris, 2016.

\bibitem{Egert-Kosmala}
M.~Egert and B.W. Kosmala, \emph{A note on complex interpolation of
  quasi-{B}anach function spaces}, In preparation.

\bibitem{Kalton-Peck-Roberts}
N.~J. Kalton, N.~T. Peck, and James~W. Roberts, \emph{An {$F$}-space sampler},
  London Mathematical Society Lecture Note Series, vol.~89, Cambridge
  University Press, Cambridge, 1984. \MR{808777}

\end{thebibliography}

\printindex

\end{document}